
\documentstyle{amsppt}
\topmatter

\parskip=0.5\baselineskip
\baselineskip=1.1\baselineskip
\parindent=0pt
\loadmsbm
\UseAMSsymbols
\input epsf
\hoffset=0.75truein
\voffset=0.5truein

\def\v{\hskip -3.5pt }
\def\epsgram|#1|#2|#3|{
        {
        \smallskip
        \hbox to \hsize
        {\hfill
        \vbox{ \vskip 6pt \centerline{\it Diagram #2}
            \centerline{\epsfbox{#3}}
        } \v \hfill } \smallskip}}

\def\epsg|#1|{
        {
        \smallskip
        \hbox to \hsize
        {\hfill
        \vbox{ \vskip 6pt
            \centerline{\epsfbox{#1}}
        } \v \hfill } \smallskip}}


\def\perim{\operatorname{perim}}
\def\area{\operatorname{area}}
\def\RR{{\Bbb R}}
\def\HH{{\Cal H}}
\def\ZZ{{\Bbb Z}}
\def\op{\operatorname{op}}
\def\arc{\operatorname{arc}}
\def\fourth{{\root 4 \of 12}}

\title The honeycomb conjecture
\endtitle
\author Thomas C. Hales \endauthor

\abstract  This article gives a proof of the classical honeycomb
conjecture: any partition of the plane into regions of equal area
has perimeter at least that of the regular hexagonal honeycomb tiling.
\endabstract

\endtopmatter

\document
\footnote""{\line{\hfill\it version - 4/17/00}}

\epsg|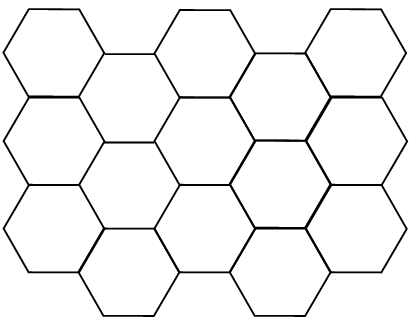|

\head 1. Introduction\endhead

Around 36 B.C., Marcus Terentius Varro, in his book
on agriculture, wrote about the hexagonal form of the bee's honeycomb.
There were two competing theories of the hexagonal
structure.  One theory held that the hexagons
better accommodated the bee's six feet.  The other theory, supported
by the mathematicians of the day, was that the structure was
explained by an isoperimetric property of the hexagonal
honeycomb.  Varro wrote, ``Does
not the chamber in the comb have six angles $\ldots$  The geometricians
prove that this hexagon inscribed in a circular figure encloses
the greatest amount of space.''

The origin of this problem is somewhat obscure.
Varro was aware of it long before Pappus
of Alexandria, who mentions it in his fifth
book.  Much of Book V follows Zenodorus's much earlier work {\it Isometric
Figures} (ca 180 B.C.).  But only fragments of Zenodorus's book remain,
and it is not known whether the honeycomb is discussed there.

The argument in Pappus is incomplete.
In fact it involves nothing more than a comparison of three suggestive cases.
It was known to the Pythagoreans that
only three regular polygons tile the plane: the triangle, the
square, and the hexagon.  Pappus states that if the same quantity of
material is used for the constructions of these figures, it is the
hexagon that will be able to hold more honey.
Pappus's reason for restricting to the three regular polygons
that tile are not mathematical
(bees avoid dissimilar figures).  He also excludes gaps between the
cells of the honeycomb without mathematical argument.  If the
cells are not contiguous ``foreign matter could enter the interstices
between them and so defile the purity of their produce'' [He,p.390].

In part because of the isoperimetric property of the honeycomb, there is
is a vast literature through the centuries mentioning the bee as
a geometer.
Thompson devotes nearly 20 pages to the literature on the
bee's cell in [Th].
Some background also appears in [K,Intro.], [Wi], [W52].
During the 18th century, the mathematical architecture of the
honeycomb was viewed as evidence of a great teleological tendency of the
universe.  Darwin explained the same
structures by natural selection. ``That motive power of the process
of natural selection having been economy of wax; that individual
swarm that wasted least honey in the secretion of wax,
having succeeded best'' [D,p.235].

The honeycomb problem has never been solved, except under special
hypotheses.
An unsolved special case of the problem is attributed to Steinhaus in [CFG,C15].
Extending Pappus's results, in 1943,
L. Fejes T\'oth proved the honeycomb
conjecture under the hypothesis that the cells are convex [FT43].
L. Fejes T\'oth predicted that a proof of the honeycomb
conjecture without the convexity hypothesis would ``involve considerable
difficulties'' [FT64a,p.183].
Elsewhere, he writes about his proof for convex cells,
``There is no doubt that the same is true for general cells.  Nevertheless,
this conjecture resisted all attempts at proving it'' [FT64b].
This paper gives a proof without the assumption of convexity.

Convexity is a highly restrictive hypothesis.
This hypothesis immediately forces the boundaries of the cells to be
polygons.
By the isoperimetric inequality,
we expect potential counterexamples
to be regions bounded by circular arcs.  One of the two regions
bounded by a positively curved arc will not be convex.  Thus, the assumption
of convexity eliminates at once almost all the candidates that should be
studied the most closely.

The geometrical properties of the three-dimensional honeycomb
cells have also been studied extensively.  The three-dimensional honeycomb
cell is a hexagonal prism built on a base of three congruent rhombuses.
The shape of the rhombic base of the three-dimensional cell
suggested the rhombic dodecahedron to
Kepler, the Voronoi cell of the face-centered cubic lattice.
During the 18th century,
many mathematicians studied the isoperimetric properties
of the base of the cells.
C. MacLaurin, in his analysis of the honeycomb,
wrote in 1743, ``The sagacity of the bees in making their
cells of an hexagonal form, has been admired of old.''
``The cells, by being hexagonal, are the most
capacious, in proportion to their surface, of any regular figures
that leave no interstices between them, and at the same time
admit of the most perfect bases'' [Mac].
In a reversal of MacLaurin's conclusions and
upsetting the prevailing opinion,
L. Fejes T\'oth discovered that the three-dimensional
honeycomb cell is not the most economical (that is, it is not
surface area minimizing) [FT64b].

The honeycomb conjecture is the two-dimensional version of the
three-dimensional Kelvin problem.  The Kelvin problem asks for the surface
minimizing partition of space into cells of equal volume.
According to Lhuilier's memoir of 1781, the problem has been described
as one of the most difficult in geometry [L,p.281].   The
solution proposed by Kelvin is a natural generalization of the hexagonal
honeycomb in two dimensions.  Take the Voronoi cells of the dual lattice
of the lattice giving the densest sphere packing.  In two dimensions, this
is the honeycomb arrangement.  In three dimensions, this gives truncated
octahedra, the Voronoi cells of the body-centered cubic.
A small deformation of the faces produces a minimal surface, which is Kelvin's
proposed solution.

Phelan and Weaire produced a remarkable counterexample
to the Kelvin conjecture.
As a result, the honeycomb problem
in two dimensions has come under increased scrutiny,
 and the need for a solution has become more acute.
F. Morgan remarks,
``In 1994, D. Weaire and R. Phelan improved on Lord
Kelvin's candidate for the least-area way to partition space into
regions of unit volume.
Contrary to popular belief, even the planar question remains open''
[M99].

It seems that
the honeycomb is minimal with respect to various optimization problems.
Even the classical problem can be expressed as a minimization of perimeter
for fixed areas or as a maximization of areas for fixed perimeters.
The first of the two presents greater difficulties and will be treated
here.  (See [FT64a,Sec. 26]).

Morgan states several versions and points out that the versions are
not known to be equivalent.  Here the situation is similar to
the sphere packing problem, which also has several competing
versions (find the densest, the solid or tight packings [CS],
the finitely stable
or the uniformly stable sphere packings [BBC]).  (Also see [Ku].)
These notions have
isoperimetric analogues.  But
here the situation is even more perplexing
because there is no upper bound on the
diameter of the cells of a partition of the plane into equal areas.
In this paper, we follow the first approach outlined by Morgan (Section 2
of his paper).  A topic for future research might be to determine to
what extent the methods of this paper can be adapted to the other
optimization problems.

Steiner's proofs of the isoperimetric problem were criticized by
Weierstrass because they did not prove the existence of a solution.
Today, general theorems assuring the existence and
regularity of solutions to
isoperimetric problems are available.  (See [T],[A],[M94], compare [B].)
This paper
depends on these results, assuring the existence of a solution to
our isoperimetric problems.

To solve the problem, we replace the planar cluster with a cluster
on a flat torus.  The torus has the advantages of compactness and
a vanishing Euler characteristic.  This part of the proof is
reminiscent of [FT43], which transports the planar cluster to a
sphere.  The key inequality, called the hexagonal isoperimetric inequality,
appears in Theorem 4.  It asserts that a certain functional is
uniquely minimized by a regular hexagon of area 1.  The isoperimetric
properties of the functional force the minimizing figure to be convex.
A penalty term prevents the solution from becoming too ``round.''
The optimality of the hexagonal honeycomb results.

\subhead Acknowledgements\endsubhead
I thank J. Sullivan for many helpful comments.
I am particularly grateful to F. Morgan for many comments, suggestions,
and corrections.
I thank D. Weaire for his recommendation,
``Given its celebrated history, it seems worth a try $\ldots$''

\head 2.  Statement of the Theorem.
\endhead

We follow [M99] in the formulation of Theorem 1-A.   Let $\pi_N = N\tan(\pi/N)$
be the isoperimetric constant for a regular $N$-gon.  That is, $4\pi_N$
is the ratio of the circumference squared to the area of a regular $N$-gon.
The particularly important case, the perimeter
$2\sqrt{\pi_6}=2\fourth$
of a regular hexagon with unit area, will be used frequently.
Let $B(0,r)$ be a disk of radius $r$ at the origin.

\proclaim{Theorem 1-A (Honeycomb conjecture)}
Let $\Gamma$ be a locally finite graph in $\RR^2$, consisting of
smooth curves, and such that
$\RR^2 \setminus \Gamma$ has infinitely many bounded connected components,
all of unit area.
Let $C$ be the union of these bounded components.
Then
$$\limsup_{r\mapsto\infty} {\perim(C\cap B(0,r))\over \area(C\cap B(0,r))}
    \ge \fourth.$$
Equality is attained for the regular hexagonal tile.
\endproclaim

The limit is insensitive to compact alterations. Therefore, there is no
uniqueness statement for the theorem in this form.  The uniqueness
of the hexagonal tile appears
in the compact version of Theorem 3 below.

Theorem 1-A has stronger hypotheses than necessary.
It assumes that the curves are piecewise smooth.  Each cell must
be connected with unit area.  There can be no {\it interstices\/} between
the cells.  Why are disks $B(0,r)$ used for the truncation?  And why
must the inequality involve $\limsup$?

We present a second version (1-B)
of the theorem that has weaker hypotheses. Before stating the theorem,
we discuss the form such a theorem might take.
Let $T_1,\ldots,T_k,\ldots$ be a countable sequence of disjoint
subsets of $\RR^2$, representing the cells of a general cluster.
It is natural to assume that
for each $i$, the topological
boundary of $T_i$ has finite $1$-dimensional Hausdorff measure.
By a result of Federer, this implies that $T_i$ is measurable and
that the current boundary $\partial T_i$ is rectifiable (see [M94,2.1],
[F,4.5.12,2.10.6]).
In general, the $1$-dimensional Hausdorff measure of $\cup_i \partial T_i$
will be infinite.  To get a finite perimeter, we truncate by fixing
a compact set $K\subset {\RR^2}$ (for example, a disk of radius
$\rho$).  Let $R_i \subset T_i\cap K$ be such that
    $\HH^1(\partial R_i\setminus \partial T_i) = 0$, where $\HH^1$ is
the $1$-dimensional Hausdorff measure.  For example, we
could take $R_i$ to be the union of connected components of $T_i$
contained in $K$.  We can measure the characteristics of the candidate
$\{T_i\}$ through the Hausdorff measure
of the sets $\cup\partial R_i$, and by
taking a compact exhaustion of the plane with sets $K$.
To state the honeycomb inequality, it is not necessary to refer to
the original cells $T_i$; it can be formulated in terms of $R_i$ and
$K$, where we now allow the area of $R_i$ to be less than $1$.
This provides motivation for the
honeycomb problem in the following general form.

\proclaim{Theorem 1-B (Honeycomb conjecture for disconnected
regions)} Let $K$ be a compact set in the plane containing
disjoint measurable sets $R_1,R_2,\ldots$. Assume that each $R_i$
has a rectifiable current boundary $\partial R_i$. Set $\alpha_i
=\min(1,\area(R_i))$.  Set $\Gamma = \cup_i \partial R_i$. Assume
    $\alpha_i>0$ for some $i$.
    Then
$$\HH^1(\Gamma) > \fourth\sum \alpha_i.$$
\endproclaim

Asymptotically, this inequality is sharp.
For example,
take $T_i$ to be the regions of the honeycomb tile, and $K=B(0,\rho)$,
a disk of radius $\rho$.  Take $R_i = T_i$, if $T_i\subset B(0,\rho)$, and
$R_i=\emptyset$, otherwise. Then $\alpha_i=0$ or $1$, and
$\sum\alpha_i$ is the number (or area) of the
hexagonal tiles contained entirely
in $B(0,\rho)$.  Moreover,
$\HH^1(\Gamma)$ is asymptotic to $\fourth\sum\alpha_i$, as
$\rho\mapsto\infty$.

We combine regions $R_i$ if the sum of their areas is
less than $1$.  This does not change $\HH^1(\Gamma)$ or $\sum\alpha_i$,
so that the validity of the theorem is not affected by combining regions.
Regions can be combined as long as there are at least two regions
with area between $0$ and $1/2$.  This allows us to assume without
generality, when we prove 1-B,
 that there are only finitely many $R_1,\ldots,R_n$.

To apply regularity results, we consider an optimization problem
that will lead to a lower bound on $\HH^1(\Gamma)$.
We fix the constants $\alpha_i\in[0,1]$, for $i=1,\ldots,n$.
Consider the
optimization problem of varying $R_1,\ldots,R_n$ so as to minimize
 $$\HH^1(\cup \partial R_i)$$
subject to the constraint that $\area(R_i)\ge\alpha_i$ for $i=1,\ldots,n$.
(We no longer require that $R_i\subset K$.)

By the existence and regularity results of [M94], there exists
$\Gamma=\Gamma(\alpha_1,\ldots,\alpha_n)$
that minimizes the 1-dimensional Hausdorff measure of
the boundary $\cup\partial R_i$ of the corresponding regions $R_i$.
The boundary $\Gamma$
consists of finitely many arcs of circles (possibly reducing to straight
lines) meeting at vertices of degree $3$.  (Morgan formulates the
optimization problem with equality
constraints $\area(R_i)=\beta_i$, where
$\beta_i\ge\alpha_i$ are fixed areas.  To get the existence and regularity
of $\Gamma(\alpha_1,\ldots)$ from this, we apply Morgan's optimization
to a set of constants $\beta_i\ge\alpha_i$ giving the shortest perimeter.)
The minimizing set $\Gamma$ is connected [C].
Each connected component
of each $R_i$ is simply connected.
Let $R_0$ be the union of the connected components of
${\RR^2}\setminus \Gamma$,
other than the components of $R_1,\ldots,R_n$.
$R_0$ is connected (otherwise remove edges between $R_i$ and $R_0$ to
shorten $\Gamma$ and increase $\area(R_i)$).

\bigskip

\proclaim{Theorem 2 (Honeycomb conjecture, finite version)}
Assume $0< \alpha_i\le 1$ for $i=1,\ldots,n$.
Let $A=\sum_i\alpha_i$.
Then
$$\perim(\Gamma(\alpha_1,\ldots,\alpha_n))
> A\fourth.$$
\endproclaim

{\bf Remark 2.1.} By an argument in [M99], the case
$A=\sum_i\alpha_i \le 398$ is elementary.
Assume we have $R_1,\ldots,R_n$ with $\area(R_i) \ge\alpha_i\in[0,1]$.
Applying the isoperimetric inequality to each $R_i$ and then applying
it once again to the union of the $R_i$, we find
$$
\align
2\perim(\Gamma(\alpha_1,\ldots,\alpha_n))
    &\ge \sum 2\sqrt{\pi\alpha_i} + 2 \sqrt{\pi A}\\
        & \ge 2( \sqrt{\pi} + \sqrt{\pi}/\sqrt{A})A \\
        & > 2\fourth A.
\endalign
$$

\proclaim{Lemma 2.2 (F. Morgan)}  Theorem 2 implies Theorem 1-A.
\endproclaim

\demo{Proof}
We consider the particular case $\alpha_1=\cdots=\alpha_n=1$.
Let
$$\rho_n= \perim(\Gamma(1,\ldots,1))/n.$$
Let $\rho_\infty$ be the infimum of the left-hand side of the inequality
of Theorem 1, as $C$ runs over all planar clusters satisfying the
conditions of the theorem.  Take $C$ be the hexagonal honeycomb tiling
to see that $\fourth\ge\rho_\infty$.  Theorem 2 gives
$$\liminf_n \rho_n\ge\fourth\ge \rho_\infty.$$
By [M99,2.1], we have $\rho_\infty \ge \liminf_n \rho_n$.  The result follows.
\qed
\enddemo

\smallskip
{\bf Remark 2.3.}  In the third edition of [M95],
Frank Morgan extends his truncation lemma [M99,2.1] to areas less
than $1$.  This permits a generalization of Theorem 1-A to
cells that are not connected.
\smallskip

\proclaim{Lemma 2.4}  Theorem 2 implies Theorem 1-B.
\endproclaim

\demo{Proof}
If $\Gamma$ is the current boundary in Theorem 1-B, then the
optimization problem described above yields
$\Gamma(\alpha_1,\ldots,\alpha_n)$. Its perimeter gives a lower bound
on the 1-dimensional Hausdorff measure of $\Gamma$.
\qed
\enddemo

{\bf Remark 2.5.}  Assume that $\Gamma$ is a finite connected
collection of analytic
arcs, and that $\Gamma$ is the
boundary of (finite unions of) simply connected bounded
regions $R_i$.
Let $\alpha_i=\min(1,\area(R_i))$, and set $A=\sum_i\alpha_i$.
Set
$$F(\Gamma,A) = \perim(\Gamma) - A\fourth.$$
Theorem 2 is false iff $F(\Gamma,A)\le0$ for some $\Gamma$, $A$.
To check Theorem 2, we can therefore make a finite number of modifications
that decrease the value of $F(\Gamma,A)$.
\bigskip

{\bf Remark 2.6.}
We claim that for the proof of Theorem 2, we can
assume that the regions $R_i$ are connected.  Each $R_i$ is a
disjoint union of connected components $R_{ij}$.  Let
$\alpha_{ij} = \min(1,\area(R_{ij}))$, and $A'=\sum_{ij}\alpha_{ij}$.
We have $A\le A'$.  Apply Remark 2.4.
\bigskip

{\bf Remark 2.7.}
Let $a_0 = 2\pi\sqrt3/3 = 4\pi/\pi_6$. If $R$ is one of the regions,
let $M=M(R)$ be the number of sides that its boundary has (counted by
the number of different regions neighboring $R$ along analytic arcs).
We claim that we may assume that each $R$ satisfies
$$\area(R)\ge a_0/M^2.$$

In fact, let $R$ be a region of area $a$ less than $a_0/M^2$.
By the isoperimetric inequality, the perimeter of $R$ is at least
$2\sqrt{\pi a}$, so it shares at least the length $(2/M)\sqrt{\pi a}$
with some neighboring region.
Let $\Gamma'$ be the collection of analytic arcs obtained by deleting
the edges shared with this neighboring region.  By deleting these edges,
we ``pop'' the bubble $R$ and combine it with the
neighboring region, so that $\Gamma'$ bounds one fewer region than $\Gamma$.
  If $A'$ is the sum of the $\alpha_i$ with respect to
$\Gamma'$, then $A'\ge A- a$.  We claim that
$F(\Gamma,A)> F(\Gamma',A')$. In fact,
$$
\align
F(\Gamma,A) - F(\Gamma',A')
         &= \perim(\Gamma)-\perim(\Gamma') + (A'-A)\fourth\\
&\ge (2/M)\sqrt{\pi a} -a\fourth\\
&= \sqrt{a}
    \left( (2/M)\sqrt{\pi}- \sqrt{a}\fourth \right)\\
&> \sqrt{a}M^{-1}
    (2\sqrt{\pi} - \sqrt{a_0}\fourth)\\
&= 0.
\endalign
$$

Deleting edges in this way may lead to regions that are not simply
connected.  However, when an edge is deleted,
arbitrarily short edges can be added to make
the regions simply connected again.  (Translate the boundary components
of $\Gamma$ so that they are arbitrarily close to one another, and
then connect them with a short edge.)  We can choose this short edge
to be so short that this modification to $(\Gamma',A')$ still gives
a value of $F$ less than the original $F(\Gamma,A)$.
\bigskip

\head 3. Honeycombs on a torus
\endhead

Let $\RR^2/\Lambda$ be a torus of area at least 1.
Take a partition of the torus into a finite number of simply
connected regions.
Assume that the boundary consists of a finite number
of simple rectifiable curves, meeting only at endpoints.

In view of Remark 2.7, we define $a(M) = \min(2\pi\sqrt3/(3M^2),1)$.
Assume the connected, simply connected regions are
$R_1,\ldots, R_n$, and
assume that $\area(R_i)\ge a(M_i)$.

\proclaim{Theorem 3 (Honeycomb conjecture on a torus)}
    $$\perim(\cup\partial R_i) \ge \sum_{i=1}^n {\alpha_i}\fourth.$$
Equality is attained if and only if every $R_i$ is a regular hexagon
of area 1 and each $\alpha_i=1$.
\endproclaim

This will follow as an immediate consequence of an isoperimetric
inequality proved in the next section.  The existence of a honeycomb
tiling depends on the shape and size of the lattice:
$\Lambda$
must be a sublattice of a lattice formed by the
tiling by unit area regular hexagonal tiles.
In particular, the area of $\RR^2/\Lambda$ must be an integer.

\proclaim{Lemma 3.1} This Theorem implies the honeycomb conjecture
for a finite number of cells (Theorem 2).
\endproclaim

\demo{Proof}
Consider $(\Gamma(\alpha_1,\ldots,\alpha_n),A=\sum \alpha_i)$
appearing in Theorem 2.
By Remarks 2.6 and 2.7, we may modify the example so
that the regions in the example are
connected, simply connected, and such that each region $R$ satisfies
the inequality
$$\area(R)\ge 2\pi\sqrt3/(3M^2)\ge a(M).$$
Let $(\Gamma',A')$ denote this modification. If we have $A'\le 398$, then
Theorem 2 follows from Remark 2.1.  Assume that $A'\ge 398$.
To complete the proof of the lemma, we will move the cluster
to the torus.  This will involve adding an additional edge of length
$\sqrt{1/A'}$ and an additional region of area at least $1$.
Theorem 3 applied to this situation will give the inequality
$$\perim(\Gamma')+\sqrt{1/A'} \ge (1+A')\fourth.$$
When $A'\ge 398$, this yields the inequality of Theorem 2.

Before moving the cluster to a torus,
we first move it to a cylinder.  Pick a diameter to the cluster
(a segment between maximally separated points $p_1$ and $p_2$
on $\Gamma$).  Then move
the cluster to the cylinder $\RR^2/\ZZ v$,
where $v$ is the translation along the
length of the diameter.  The map of the cluster to the cylinder is
injective, except at the points $p_1$ and $p_2$, which become identified.
Since $p_1$ and $p_2$ are maximally separated, the cluster fits inside
a square of edge length $|v|$ with a pair of sides parallel to $v$.
Thus the area of the cluster is at most the area $|v|^2$ of the square.
To simplify notation, we now drop the primes from $(\Gamma',A')$.
This gives $|v|\ge \sqrt{A}$, with $A=\sum\alpha_i$, and
$\alpha_i = \min(1,\area(R_i))$.

Let $w$ be a unit vector perpendicular to $v$.  Pick $\mu>0$
to be the largest real number for which
$\Gamma+\mu w$ touches $\Gamma$ without overlap.
  Let
$$\Lambda = \ZZ v + \ZZ (\mu +\sqrt{1/A}) w.$$
The cluster descends to $\RR^2/\Lambda$, injective except at $p_1$ and $p_2$.
We add a segment of length $\sqrt{1/A}$ to join the cluster with
its translate.
On the torus, the region ``at infinity'' becomes
simply connected.  Call it $R_0$.
  By adding this extra edge, we
avoid the complications of a component with a loop representing
a nontrivial homology class in $\RR^2/\Lambda$.

Since $|v|\ge\sqrt{A}$, and the component $R_0$ has height at
least $\sqrt{1/A}$ in the direction of $w$ at every point, the
area of $R_0$ is at least 1.
Now we have a partition of the torus $\RR^2/\Lambda$ into
connected, simply connected regions of total area at least $1+A$.
\qed
\enddemo

{\bf 3.2. Torus modifications.\/}
The combinatorial structure is described
by a finite torus graph in which each face is simply connected.

Loops (edges joined at both ends to the same vertex) can be eliminated
from the graph as follows.  If the vertex has degree $>3$, then
it can be considered a limit of multiple vertices of degree
$3$ and edges of length $0$.  This can be done in such a way that
the vertex on the loop has multiplicity at least $2$.
If the vertex has degree $3$, let $e$ be the other edge that meets the
loop at the vertex.  Both sides of the edge $e$ bound the same region $E$.
 Removing $e$ leads to a
non-simply connected component.  The loop $\Gamma_1$
can be moved arbitrarily
close to another boundary component $\Gamma_2$.
It can be joined to $\Gamma_2$
with two edges (one edge moving from a point vertex $v_1$
of $\Gamma_1$ down to $w$ on $\Gamma_2$, and another edge moving back to
$v_2$ in the opposite direction, as in Diagram 3-a).  This decreases
the area of $E$, which is not permissible.  But by scaling the entire
torus by a homothety ($x+\Lambda\mapsto tx+t\Lambda$),
its area is restored.  Since we can make the decrease
in area arbitrarily small, we can make the increase in perimeter due to
the homothety arbitrarily small.  In particular, we can arrange the
the total increase in perimeter by this process is offset by the
length of the edge $e$ that was removed.

\smallskip
\epsgram|2|3-a|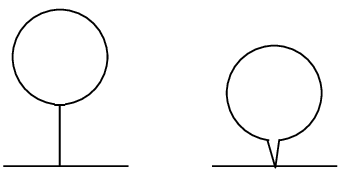|
\smallskip

(Alternatively, we could let $v_1=v_2$, and let the return path to $v_1$
be the same as the path from $v_1$ to $w$ in the opposite direction.
Then $v_1$ is to be considered a vertex of multiplicity $2$, etc.)

If the torus graph has any vertices of degree greater than $3$, we view
them as vertices of degree $3$ with degenerate edges of length $0$.
We can do this in a way that avoids creating any new loops.
Each vertex has degree $3$.

We are now in the situation where the boundary of every region $P$ is a
combinatorial $N$-gon, for some $N\ge2$, where $N=N(P)$ is the number of
directed edges bounding the region.
(The boundary of $P$ might traverse a segment twice in opposite directions.)
We have the Euler characteristic of the torus graph
$$0 = V-E+F = \sum_P (1-N(P)/6),\tag 3.3$$
where $P$ runs over the regions.

\bigskip

\head 4. A hexagonal isoperimetric inequality for closed plane curves.
\endhead

Let $\Gamma$ be a closed
piecewise simple rectifiable curve
 in the plane.
 In our application, we will take
$\Gamma$ to be a lift from the
torus to the plane of one of the $\partial R_j$ from Section 3.

We use the parametrization of the curve to give it a direction,
and use the direction to assign a signed area to the bounded
components of the plane determined
by the curve.
For example, if $\Gamma$ is a piecewise
smooth curve, the signed area is
given by Green's formula
$$\int_{\Gamma} x dy.$$
Generally, we view $\Gamma$ as an integral current [M95,p.44].
We let $P$ be an integral current with boundary $\Gamma$.
(In applications, $P=R_j$, for some $j$.)
Expressed differently,
we give a signed area by assigning an multiplicity $m(U)\in\ZZ$
 to each bounded
component $U$ of $\RR^2\setminus\Gamma$. (An illustration appears
in Diagram 4-a.)
The area is $\sum m(U)\area(U)$.
$P$ is represented by the formal sum $P = \sum m(U) U$.

Let $v_1,\ldots,v_t$, $t\ge 2$, be a finite list of points on $\Gamma$.
We do not assume that
the points are distinct.
Index the points $v_{1},v_{2},\ldots,v_{t}$, in the order
provided by the parametrization of $\Gamma_i$.  Join $v_{i}$
to $v_{i+1}$ by a directed line segment $f_{i}$
(take $v_{t+1}=v_{1}$).
The chords $f_i$ form a generalized polygon, and from the direction
assigned to the edges, it has a signed area $A_P\in\RR$.

\smallskip
\epsgram|2|4-a|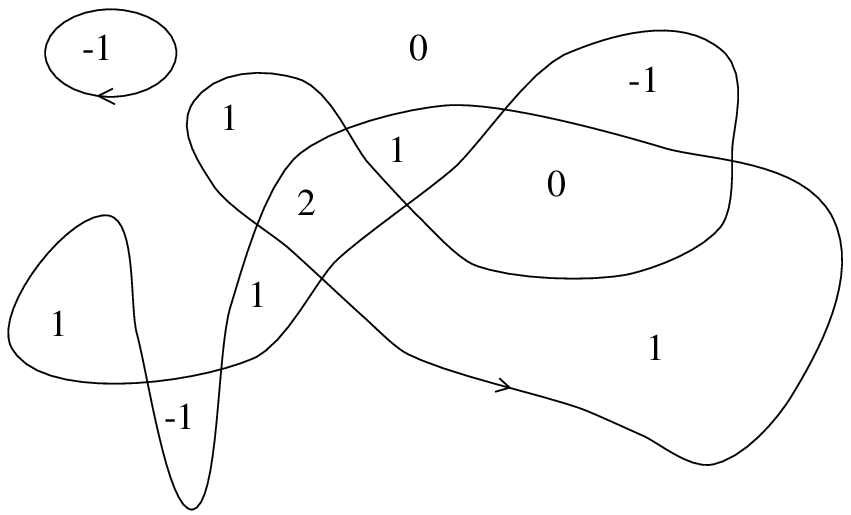|
\smallskip

Let $e_i$ be segment of $\Gamma$ between $v_{i}$
to $v_{i+1}$.  Let $f^{\op}$ be the chord $f$ with the
orientation reversed.
Let $x(e_{i})\in\RR$ be the signed area of the integral current
bounded by $(e_{i},f_{i}^{\op})$.
Let $E(P)=\{e_{i}\}$ denote the set of edges of $P$.
Accounting
for multiplicities and orientations, we have
$$\area(P) = A_P + \sum_{e\in E(P)} x(e).$$
Let $\alpha(P) =\min(1,\area(P))$.

Define a truncation function $\tau:\RR\to \RR$ by
$$\tau(x) = \cases
        1/2, & x\ge 1/2,\\
        x,   & |x|\le 1/2,\\
        -1/2, &x\le -1/2.\endcases
$$
Set $\tau_0 = 1/2$.  Set $T(P)= \sum_{E(P)}\tau(x(e))$.
Recall that the perimeter of a regular hexagon of unit area
is $2\fourth$.
Let $L(P)$ be the length of $\Gamma$.
Let $N(P)$ be the number of points $v_i$ on $\Gamma$, counted
with multiplicities.
Recall from Remark 2.7 that $a(N) = \min(2\pi\sqrt3/(3N^2),1)$.

\proclaim{Theorem 4 (hexagonal isoperimetric inequality)}
Define $P$, $L(P)$, $N(P)$, and $a(N)$ as above.
Assume that the signed area of $P$ is at least $a(N(P))$.
 Then
$$L(P)\ge - T(P)\fourth - (N(P)-6) 0.0505 + 2\alpha(P)\fourth.$$
Equality is attained if and only if $P$ is a regular hexagon
of area 1.
\endproclaim

The theorem will be proved below.

\proclaim{Lemma 4.1}  The hexagonal
isoperimetric inequality implies the honeycomb
conjecture for a torus (Theorem 3).
\endproclaim

\demo{Proof}  We will apply this inequality to the $n$ different
regions $R_1,\ldots, R_n$ in the torus partition.
We let the points $v_i$ be the endpoints of the simple
curves described in Section 3.

Let $P$ be one of the regions $R_i$.
The regions in Theorem 3 satisfy $\area(P)\ge a(M(P))$, where $M$
is the number of different regions bounding $P$.  We have $N\ge M$,
so $\area(P)\ge a(M(P)) \ge a(N(P))$.  Thus, $P$ satisfies the
area constraint of the hexagonal isoperimetric inequality.

Each edge of $P$ occurs
with opposite orientation $e^{\op}$ on a neighboring region $P'$,
and by construction $x(e) + x(e^{\op}) = 0$.
Thus, summing over all directed
edges of the partition, we get
$$\sum_e x(e) = 0.$$
Since $\tau$ is an odd function,
we have $\tau(x(e))+\tau(x(e^{\op}))=0$, so also
$$\sum_e \tau(x(e)) = \sum_P T(P) = 0.\tag{4.2}$$

By Theorem 4,
$$
\align
2\perim&(\Gamma)=
    \\
    \sum_P &L(P) \ge -\fourth\sum_P T(P)
    + 6 (0.0505) \sum_P (1-N(P)/6)
    +2\sum_P \alpha(P)\fourth.
\endalign
$$
Thus, the lemma follows from Euler (Equation 3.3) and Equation 4.2.
\qed\enddemo

Set
$$
\align
\Delta(P) &= L(P) + \epsilon(N(P),\alpha(P),T(P)),\\
\epsilon(N,\alpha,T) &= T \fourth + (N-6)0.0505 - 2\alpha\fourth,\\
X(P) &= \sum_{E(P)} x(e).
\endalign
$$

{\bf Remark 1.}
The inequality is false without the truncation.  For example,
let $P$ be constructed as a simple closed curve
of area 1 bounded by 3
inverted circular arcs of the same curvature and the same length
(Diagram 4-b).
When the circular arcs are sufficiently long,
$L(P) + \epsilon(3,1,X(P)) <0$.

\smallskip
\epsgram|2|4-b|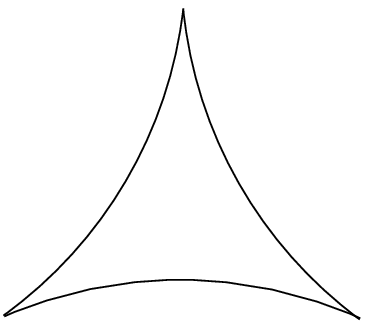|
\smallskip

{\bf Remark 2.}
If the region $P$ is a polygon of area 1 with $v_i$ as vertices,
then $\alpha(P)=1$, $x(e) = 0$ for all $e$, and $T(P)=0$.
The inequality in this case is essentially the one used
by Phelan and Weaire in [Ph].

{\bf Remark 3.}
The sharp case of the inequality occurs for $N(P)=6$, $\alpha(P)=1$,
 and $x(e)$ near $0$, so that $X(P) = T(P)$.
$y= -\epsilon(6,1,X)$ is the tangent line to $y=L(P_X)$ at $X=0$
of the following region $P_X$ (Diagram 4-c).
Take a regular hexagon of area
$1-X$ and add six circular arcs of the same curvature to the
six edges to make the total area 1 (inverting the arcs if $X<0$).

\smallskip
\epsgram|2|4-c|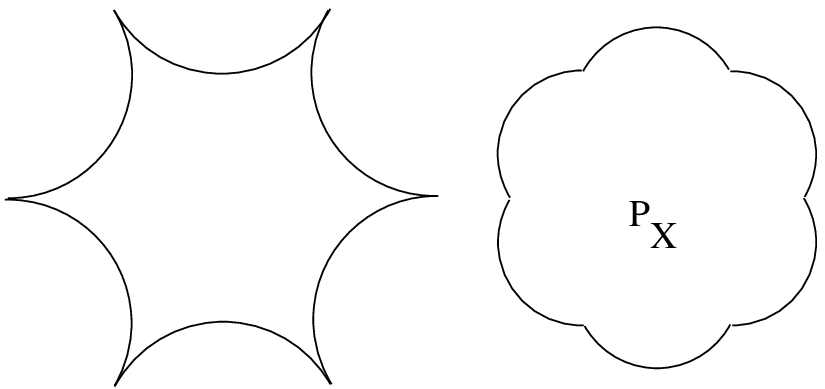|
\epsgram|2|4-d. The graph of $y=L(P_X)$ and its tangent line $y=-\epsilon(6,1,X)$|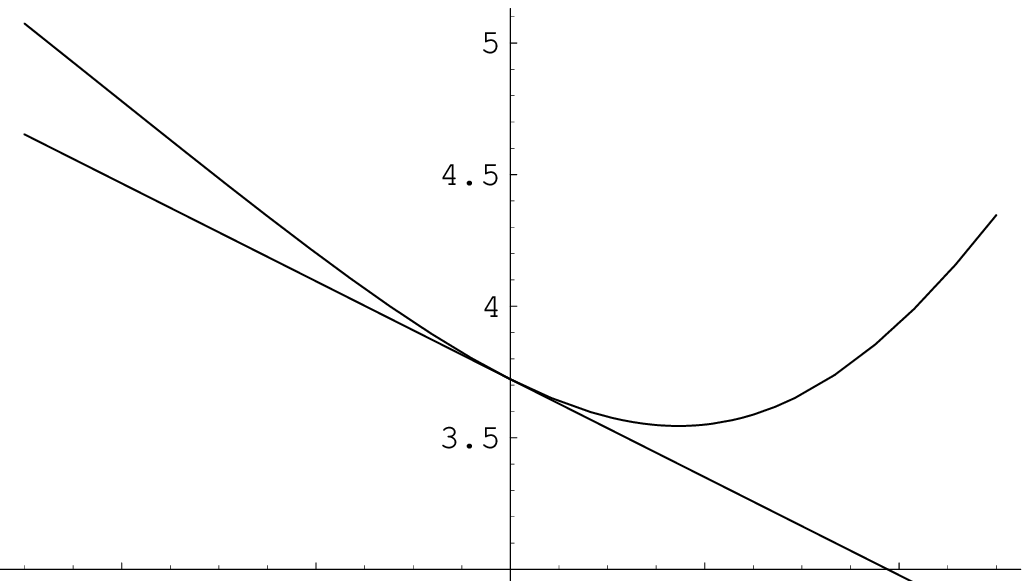|
\smallskip

\head 5. Preliminaries for the Proof
\endhead

Without loss of generality, we modify $P$ to decrease the
perimeter, maintaining the lower bound on the area, and holding
$\epsilon$ fixed.  Or we may also modify $P$ by fixing the
perimeter, maintaining the lower bound on the area, and
decreasing $\epsilon$.

For each chord $f$, we may apply the isoperimetric inequality
to the curve $(e,f^{\op})$ to replace $e$ with a circular
arc with the same enclosed signed area.  The isoperimetric inequality
for integral currents appears in [F;4.5.14].  Uniqueness and regularity
follow along the lines of [M94].

We may replace the polygon $(f_1,\ldots,f_N)$ by a convex polygon
with the same edge lengths in the same order, that has at
least the area as the original polygon.  Thus, we may assume that
the polygon has positive area.

\smallskip
\epsgram|2|5-b. A deformation reducing the perimeter at fixed area.|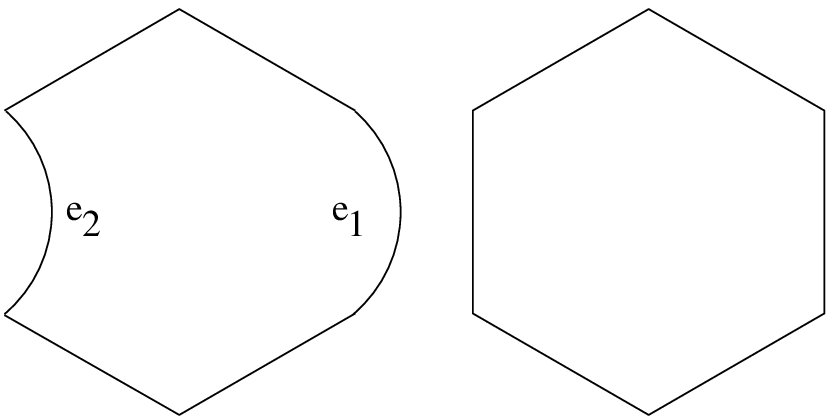|
\smallskip

If any $x(e)<-\tau_0$, we decrease the curvature of $e$ until $x(e)=-\tau_0$.
This leaves $T(P)$ unchanged, decreases the perimeter, and increases
$\area(P)$.  Thus, there is no loss in generality, if we assume
$x(e)\ge -\tau_0$ for all $e$.

We leave edges $e$ with $x(e)>\tau_0$ untouched.

If there are two edges $e_1,e_2$ of $P$
with $0<x(e_1)\le\tau_0$
and $x(e_2)<0$ (and so also $x(e_2)\ge-\tau_0$),
we deform
$P$ by decreasing the curvature of the arcs $e_1$ and $e_2$
preserving $x(e_1)+x(e_2)$, thereby decreasing the perimeter
$L(P)$ (see Diagram 5-b).
$T(P)$ is unchanged.  Continuing in this way, we may assume without
loss of generality that all $x(e)\in[-\tau_0,\tau_0]$ have
the same sign.  We consider two cases depending on whether
there is an edge $e$ with $x(e)>\tau_0$.

(Case I) For all $e$, $x(e) \ge -\tau_0$, and there exists $e$ such
that $x(e) > \tau_0$.  All $x(e)\in[-\tau_0,\tau_0]$ have the same sign.

(Case II) For all $e$, $|x(e)|\le\tau_0$.  All $x(e)$ have the same
sign.

The theorem will be proved by a separate argument for the two cases.

\head 6. Lower bounds on $L(P)$.
\endhead

We were not able to find a single estimate of $L(P)$ that leads
to the Theorem in all cases.  Instead, we rely on number of lower
bounds on $L(P)$ (and $\Delta(P)$).  Most are based on
the isoperimetric inequality.

\subhead Lower bound $L(N,\alpha,X),$ $N=3,4,\ldots$
\endsubhead
The perimeter of a regular
$N$-gon of area 1 is
$2\sqrt{\pi_N}$, where $\pi_N = N\tan(\pi/N)$.
 The polygon $(f_1,\ldots,f_N)$ has area
at least $\alpha-X$.  By the isoperimetric inequality for polygons, it has perimeter
at least
    $$L(N,\alpha,X) := 2 \sqrt{(\alpha-X)\pi_N}.$$
Each arc $e_i$ has length at least that of $f_i$, so $L(P)\ge L(N,\alpha,X(P))$.

\subhead Lower bound $L_+$
\endsubhead
By the isoperimetric inequality,
$$L(P)\ge L_+(\alpha):= 2\sqrt{\pi\alpha},$$
the perimeter of a circle of area $\alpha$.

\subhead Lower bound $L_-$
\endsubhead
Let $J$ be the set of indices of the edges with $x(e_j)<0$, for $j\in J$.
Let $X_J = \sum_J x(e_j)$.
By reflecting each edge $e_j$, $j\in J$, across the chord $f_j$,
$P$ is replaced with a region of the same perimeter and area
at least $\alpha-2 X_J$.  By the isoperimetric inequality,
$$L(P)\ge L_-(\alpha,X_J) = 2\sqrt{(\alpha-2X_J)\pi}.$$

\subhead Lower bound $L_D$
\endsubhead
Assume here that $|x(e)|\le \tau_0$, for each $e$.
Dido found the curve of minimum length with both endpoints
on a given line, subject to the condition that the curve and the line
bound a given area.  The solution is a semicircle cut through the center
by the line.
Applying this to an arc $e$, we find that the length of $e$
is at least $\sqrt{2\pi|x(e)|}\ge |x(e)|\sqrt{2\pi/\tau_0}$.
Thus, $$L(P) \ge L_D(X_D) := X_D \sqrt{2\pi/\tau_0},$$
where $X_D =\sum_e |x(e)|$.  In Case II, $X_D = |T|=|X|$.

\subhead Lower bound $L'(N,\alpha,X)$
\endsubhead  Only this last bound is new.
For $\ell,x\ge0$,
let $\arc(\ell,x)$ be the length of a circular arc chosen so that
together with a chord of length $\ell$ joining its endpoints, the
enclosed area is $x$.  For example, $\arc(\ell,0)=\ell$ and
$\arc(0,x) = 2\sqrt{x\pi}$.
Let $L(N,\alpha,x)$ be as above.
Let
$$L'(N,\alpha,X) = L(N,\alpha,X) \arc(1,|X|/L(N,\alpha,X)).$$

\proclaim{Proposition 6.1}
If all the chords $f_i$ of $P$ have length at most 1,
if $|X|\le 0.119$, $0.996\le\alpha\le1$, and if $N\le7$, then
    $$L(P)\ge L'(N,\alpha,X),$$
 where $X=X(P)$, $N=N(P)$, $\alpha=\alpha(P)$.
\endproclaim

This will be proved in Appendix 1.

\subhead 6.2. Equal curvature condition
\endsubhead
Here is a simple observation about the lower bounds on perimeters
that we will refer to as the {\it equal curvature condition.}
(A version for polygons was known to Zenodorus [He,p.210].)

Suppose that we have two chords $f_1$ and $f_2$ of circular arcs
$e_1$ and $e_2$.  Minimize
the sum of the lengths of $e_1$ and $e_2$, fixing $f_1$ and $f_2$, and
constrained so the
sum of the two enclosed areas is fixed.
Two arcs of equal curvature give the minimum.
If an arc is
more than a semicircle, it occurs
along the chord of greater length
(or one of the two if the chords are equal in length).

To see this result, form a triangle with the two chords and a third edge
of variable length $t$.  Adjust $t$ until the circumscribing
circle gives arcs of the correct combined enclosing area on the two
chords. Any shorter perimeter contradicts the isoperimetric
inequality.

\head 7.  Case I of the proof of the hexagonal isoperimetric inequality
\endhead

\subhead Digons
\endsubhead
Before treating Case I, we treat the case of digons separately for
both cases I and II.  Here $A_P=0$ and $N(P)=2$, so
    $$\area(P) = x(e_1)+x(e_2)\ge \alpha.$$
Also, $\alpha\ge a(N)=a(2)> 1/4$.
In Case I, $\sum x(e_i) = x_{\max} + x_{\min} \ge \tau_0-\tau_0\ge0$.
In Case II,
$T(P) = \sum \tau(x(e_i))= \sum x(e_i) \ge\alpha$.

If $T(P)>0.21$, then $\Delta(P)>0$, by the bound $L_+$:
$L(P)\ge 2\sqrt{\pi\alpha}\ge 2\alpha\sqrt{\pi}$.

Assume $T(P)\le 0.21$.  We are now in Case I, so $x_{\max}>\tau_0=\tau_{\max}$,
$$x_{\min}\le \tau_{\min} \le 0.21 - \tau_{\max} = -0.29.$$
By reflecting the arc $e$ corresponding to $x_{\min}$ across the corresponding
chord $f$, the area becomes at least $(\alpha-2x_{\min})$, without changing
the perimeter.
We then have,
$$L(P)\ge 2\sqrt{\pi(\alpha-2x_{\min})} \ge 2 \sqrt{\pi(\alpha+2(0.29))}.$$
$\Delta(P)>0$ follows.

For the rest of the proof, we assume $N(P)\ge 3$, so that in particular,
$\epsilon(N,\alpha,T)\ge \epsilon(3,\alpha,T)$.

\subhead Case I \endsubhead
We assume that for all $e$, $x(e)\ge -\tau_0$ and that for some $e$,
$x(e)>\tau_0$.  All $x(e)$ satisfying $|x(e)|\le\tau_0$ have the same sign.
In treating this case, we only need to assume that $\alpha\ge0$, rather
than $\alpha\ge a(N)$.

Assume $T(P)>0.177$.  We have the bounds
$L_+: L(P)\ge 2\alpha\sqrt{\pi}$ and
$\epsilon(N,\alpha,T)\ge \epsilon(3,\alpha,0.177)$.
It follows that $\Delta(P)>0$.

Assume next that $T(P)<-0.36$.  There
exists $x(e)<0$.
Index so that $x(e_i) >\tau_0$ for $i\in I$ and $x(e_j) \le 0$ for $j\in J$.
Set $X_I = \sum x(e_i)$, $X_J = \sum x(e_j)$, so that
$X(P) = X_I + X_J$.  The area of $P$ is
    $X_I + X_J + A_P \ge \alpha$.
Let $k=|I|$.  We have
$T(P) = k \tau_0 + X_J$.
By Dido,
$$L(P) \ge \sum_I \sqrt{2\pi|x_i|} + \sum_J \sqrt{2\pi|x_j|}
    \ge k \sqrt{2\pi \tau_0} - X_J \sqrt{2\pi/\tau_0}.$$
Then
$$\Delta(P)\ge k\sqrt{2\pi\tau_0}-X_J\sqrt{2\pi/\tau_0}
    +\epsilon(3,1,k\tau_0+X_J).$$
Substituting the upper bound $X_J\le -0.36-k\tau_0$ for $X_J$,
and then the lower bound $k\ge1$ for $k$, we find that
$\Delta(P)>0$.

Assume finally that $T(P)\in[-0.36,0.177]$.
With the same notation, we have
    $T(P) = k\tau_0 + X_J\le 0.177$, which gives $X_J\le -0.323$.
Reflecting the arcs corresponding to
negative signed areas as above, we get
$$L(P) \ge 2\sqrt{\pi(\alpha+2(0.323))}.$$
This gives $\Delta(P)>0$.

\head 8.  Case II of the proof of the hexagonal isoperimetric inequality
\endhead

Assume that $|x(e)|\le \tau_0$ for all $e$, and
that all $x(e)$ have the same sign.  Then $X(P)=T(P)$ and $X_D=|X(P)|$.

First, we will treat the case $\alpha\in[2\pi\sqrt3/(3N^2),1/4]$, and then we will
treat the case $\alpha\ge 1/4$.  If
$1/4\ge\alpha\ge 2\pi\sqrt3/(3N^2)$, we have $N\ge 4$.
If $T\ge 0$, we use $L(P)\ge L_+(\alpha)$ to get
$$
\Delta(P) \ge 2\sqrt{\pi\alpha} + (N-6)0.0505 - 2\alpha\sqrt{\pi}.
$$
The second derivative in $\alpha$ is negative, so it is enough to
check that this is positive for $\alpha = 2\pi\sqrt3/(3N^2)$, and $\alpha=1/4$.
This is elementary.

Next, if $T\in[-2.4,0]$, we use $L(P)\ge L_-(P)$.  We show that the
following is positive:
$$2 \sqrt{\pi(\alpha-2 T)} + T\fourth + (N-6) 0.0505 - 2\alpha\fourth.$$
The second derivative in $\alpha$ is negative as well as that for $T$.
Hence, it is enough to check positivity for $T=0, -2.4$,
and $\alpha=2\pi\sqrt3/(3N^2), 1/4$.  Again, the verification is elementary.

Finally, if $T\le -2.4$, we use $L(P)\ge L_D$.  It is clear that
$$\Delta(P)\ge -T\sqrt{2\pi/\tau_0} + \epsilon(N,1,T) \ge 0.$$
This completes our discussion of the case $\alpha\in[2\pi\sqrt3/(3N^2),1/4]$.

\bigskip

The rest of Section 8 is devoted to the case $\alpha\ge 1/4$.
We pick a lower bound $\tilde L(N,\alpha,X)$ for $L(P)$ from the stock of lower
bounds developed in Section 6,
according to the following schematic in the $(N,X)$ plane.

\smallskip
\epsgram|2|8-a|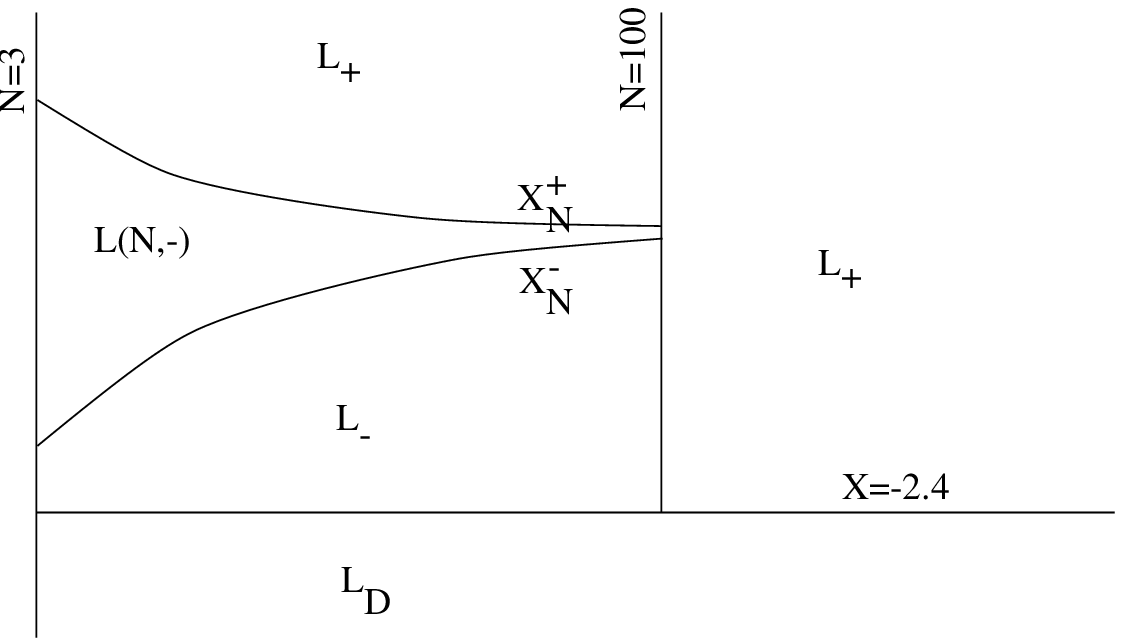|
\smallskip

The boundary between $L_+(1)$ and $L(N,\cdot)$ is the curve
$X_N^+= 1-\pi/\pi_N>0$, determined by the condition
$L_+(1)=L(N,1,X_N^+)$.
(For $N=3$, we set $X_3^+ = 0.177$, instead of $1-\pi/\pi_3$.)
Also, set $X_N^- =(-\pi+\pi_N)/(-2\pi+\pi_N)<0$,
which satisfies $L_-(1,X_N^-)=L(N,1,X_N^-)$.
When $X<0$, we have $X=X_J$, and $L_-(\alpha,X_J) = L_-(\alpha,X)$.
We omit the details of a routine calculation that shows
$$
\tilde L(N,\alpha,X)+\epsilon(N,\alpha,X)\ge \tilde L(N,1,X)
    +\epsilon(N,1,X),
$$
for $N>4$.

The function $\tilde L(N,1,X)+\epsilon(N,1,X)$
is piecewise analytic and
has a negative second derivative in $X$.
This makes it trivial to check that this function is positive
on given analytic interval by checking the values at the endpoints.
The function $\tilde L(N,\alpha,X)+\epsilon(N,\alpha,X)$, for
$N=3,4$, is also easily checked to be positive.
We find that
$$\Delta(P)\ge \tilde L(N,\alpha,X) +\epsilon(N,\alpha,X)>0,$$
except in the following two situations that will be treated below.

(1) $N=6$, $\alpha\in[0.996,1]$, and $X\in[-0.119,0.1]$,

(2) $N=7$, $\alpha\in[0.996,1]$, and $X\in [-0.082,0.0684]$.

Diagrams 8-b and 8-c show
 $\tilde L(N,1,X)+\epsilon(N,1,X)$, for $N=6,7$.

\smallskip
\epsgram|2|8-b. The graph of $\tilde L(6,1,X)+\epsilon(6,1,X)$.|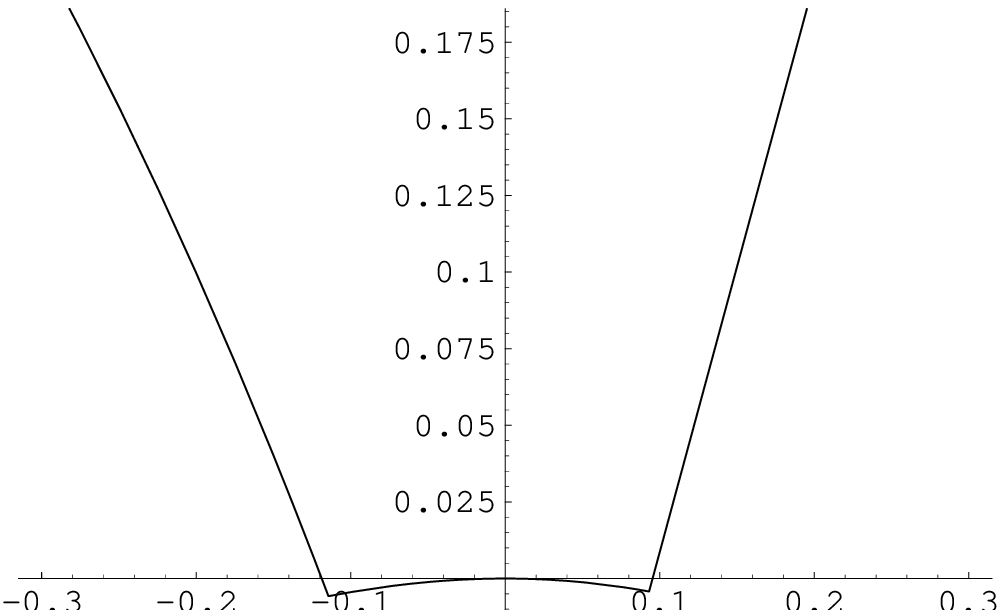|
\epsgram|2|8-c. The graph of $\tilde L(7,1,X)+\epsilon(7,1,X)$.|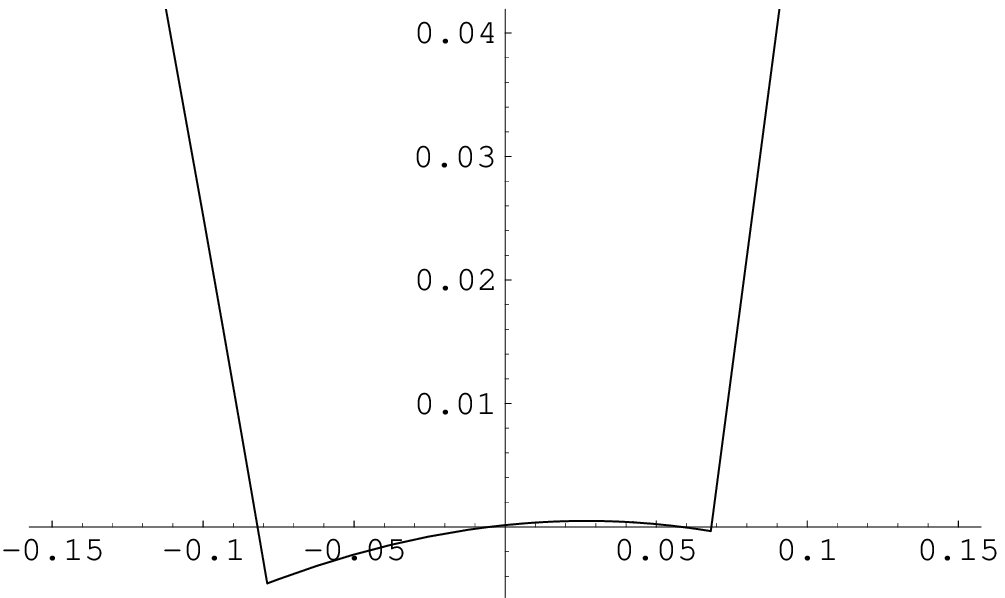|
\smallskip

For the rest of the argument, we assume we are in one of these two
situations.
If the lengths of the chords $f_i$ are at most $1$,  we
have
$$\Delta(P) \ge L'(N,\alpha,X) +\epsilon(N,\alpha,X) \ge 0.$$
The inequality on the left is Proposition 6.1.

The inequality on
the right
was checked by taking a Taylor approximation at $X=0$ (and by taking
the worst-case values for $\alpha$).  We give a few details of the
approximation in the
special case $N=6$, $|X|\le0.06$.  The constant $\alpha$ satisfies
$0<\alpha\le1$.  We also assume that $0\le\theta\le0.5$,
where $\theta$ is a variable defined below.  The function takes the form
$$L'(6,\alpha,X)+\epsilon(6,\alpha,X) =
    \fourth (2q(\theta) \sqrt{\alpha-X} - 2\alpha + X),$$
where $q(\theta) =\theta/\sin\,\theta$, and $\theta\ge0$ is defined
implicitly by the equation
$$p(\theta) := {\theta-\sin\theta\cos\theta\over\sin^2\theta} =
    {2|X|\over\sqrt{\alpha-X}\fourth}.$$
(This formula for $L'$ is presented in greater detail in the Appendix.)
We show that
    $$2q(\theta)\sqrt{\alpha-X} - 2\alpha + X\ge0,$$
with equality exactly when $X=0$.  The following inequalities are
easily verified, when the variables lie in the indicated intervals.
$$\sqrt{1-X} \ge 1 - X/2 - 0.22 X^2.$$
$${2|X|\over \fourth \sqrt{1+0.06} } \le
    {2 |X|\over \fourth \sqrt{\alpha-X}} = p(\theta) \le 0.87\theta.$$
$$q(\theta) \ge 1+{\theta^2\over6} \ge 1+
    {4 X^2\over 6(0.87^2)(1.06)\sqrt{12}}.$$
$$2q(\theta)\sqrt{\alpha-X}-2\alpha+X
    \ge 2 (1+ {4 X^2\over 6(0.87^2)(1.06)\sqrt{12} })(1-X/2-0.22X^2) - 2+X = X^2 f,$$
where $f$ is a quadratic polynomial in $X$ taking positive values for
$|X|\le 0.06$.  It follows that if the function vanishes, then $X=0$.
This implies $\theta=0$, and $q(\theta)=1$.  We then have
$$L'(6,\alpha,0)+\epsilon(6,\alpha,0) = \fourth (2\sqrt{\alpha}-2\alpha).$$
If this vanishes, then $\alpha=1$.

This completes the hardest case.  Similar calculations are left to the
reader when $N=7$, $|X|\ge 0.06$, or $\theta\ge0.5$.
The result is that
$\Delta(P)=0$, if and only if $X=0$ and
$P$ is a regular hexagon with area $1$. This is the tight case
of the hexagonal isoperimetric inequality. (The graph of $y=L'(6,1,X)$
lies between the graphs of $y=L(P_X)$ and $y=-\epsilon(6,1,X)$ in
Diagram 4-d.)

Now assume that some chord $f_i$ has length at least $1$.
The area of the polygon $P_f = (f_1,\ldots,f_N)$ is
at least $\alpha-X$.  A lower bound on $L(P)$ is the perimeter $L(P_f)$ of this
polygon.

\proclaim{Lemma}  $L(P_f) + \epsilon(N,\alpha,X) >0$, in situations (1) and (2).
\endproclaim

\demo{Proof}
We develop an isoperimetric inequality for polygons of area
at least $\alpha-X$, constrained so that one of the edges has length at least 1.
We may assume that the area is $0.996-X$.
By well-known principles,
the optimal polygon is inscribed in a circle with unconstrained
edge lengths $t$ and constrained edge length $\max(1,t)$.
The area $A_N(r)$ and perimeter $\max(1,t)+(N-1)t$ are monotonic increasing
functions of the circle's radius $r$.

Now $X=X_N(r)=0.996-A_N(r)$, so
$$\Delta(P) \ge g_N(r):=\max(1,t) + (N-1) t  + \epsilon(N,1,X_N(r)).$$
The function $g_N(r)$
is easily estimated because of the monotonicity of $t$ and $A_N$.
For $a<b$, we write
$$g_N(a,b) = \max(1,t(a))+(N-1)t(a) + \epsilon(N,1,0.996-A_N(b)).$$
Then $\Delta(P)\ge g_N(r)\ge g_N(a,b)$, for $r\in [a,b]$.
In situation (1),
    $$X_6(0.671)< -0.119\le X \le 0.1 < X_6(0.61).$$
Thus, it can be seen that $\Delta(P)>0$, by computing the constants
    $$g_6(0.61+0.001\,k,0.611+0.001\,k)>0,$$ for $i=0,\ldots,60$.
(The constants are all at least $0.02$.)
The situation (2) is similar.
\qed
\enddemo

\head Appendix 1. A proof of Proposition 6.1.
\endhead


Let $P$ have chords $f_i$ of length $\ell_i$.
In this appendix, we give a proof of the following result.

\proclaim{Proposition 6.1-A}
Choose constants $s$ and $\ell$ so that
$\ell_i\le s<\ell \le \sum \ell_i$, for all $i$.
If $|X(P)|\le \pi s^2/8$,
then
    $$L(P)\ge \ell \arc(s,|X(P)|/\ell).$$
\endproclaim

We obtain the version of the proposition that is stated
in Section 6 by taking $s=1$, $\ell=L(N,\alpha,X(P))$,
$N(P)\le 7$, and $\alpha\in[0.996,1]$.
Note that under these conditions,
$0.119<\pi s^2/8$ and $\ell \arc(s,|X(P)|/\ell) = L'(N,\alpha,X)$.

\demo{Proof}
We flatten
out the perimeter of the polygon by arranging its edges
$f_i$ along a line as shown in Diagram A-a.
Without loss of generality, we may assume that all $x(e)$ have
the same sign.  (See Section 5, noting that the truncation
$\tau_0$ is not used in this appendix.)

\smallskip
\epsgram|2|A-a|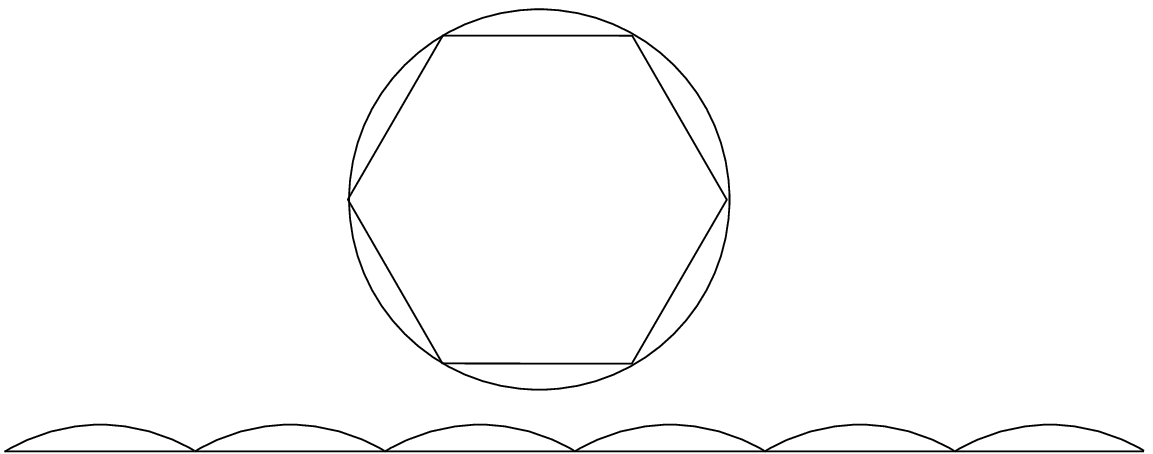|
\smallskip

Next we analyze two particular edges of lengths $(u,v)=(\ell_i,\ell_j)$.
Let the lengths of the
two circular arcs be $2\rho \theta_u$ and $2\rho \theta_v$,
where $1/\rho$ is common curvature of the the two arcs (cf. Section 6.2).
That is, $\theta_u$ is the angle subtended by the arc and its chord.

\smallskip
\epsgram|2|A-b|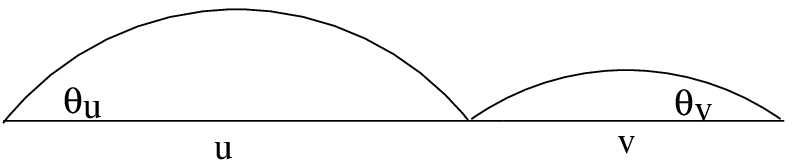|
\smallskip

Set
$$
\align
p(t) &= {t - \sin(t)\cos(t)\over \sin^2(t)},\\
q(t) &= \quad {t\over \sin(t)}.\\
\endalign
$$
We have
$$
\align
\ell_i+\ell_j &= u+v,\tag A.1\\
4(|x(e_i)|+|x(e_j)|) &= u^2 p(\theta_u) + v^2 p(\theta_v),\tag A.2\\
\perim(e_i)+\perim(e_j) &= u q(\theta_u) + v q(\theta_v),\tag A.3\\
0 &= u\sin\theta_v - v\sin\theta_u.\tag A.4
\endalign
$$
These four equations give the length of the chords, the enclosed area,
the arc length of the circular arcs, and the equal curvature condition
for the two edges.
Let $\xi(u)$ be the area $|x(e_i)|+|x(e_j)|$, viewed
as a function of $u$ by fixing the combined length $\ell_i+\ell_j$,
perimeter $\perim(e_i)+\perim(e_j)$, and the equal curvature condition.
It is defined implicitly by the equations A.1--A.4.

\proclaim{Lemma A.5}
$$\xi'(u)=\rho(\cos\theta_v-\cos\theta_u).$$
In particular, $\xi'(u)\ge0$, if $0\le\theta_v\le\theta_u\le\pi$.
\endproclaim

\demo{Proof}
We compute the derivative $\xi'(u)$ by implicit differentiation of the
equations A.1-A.4.  The differentials of
these four equations give four homogeneous linear relations
among $du$, $dv$, $d\theta_u$, $d\theta_v$, and $d\xi$.
Solving the linear system for $d\xi/du$,
we obtain the result.
\qed
\enddemo

{\bf Remark.}
Stewart Johnson and Frank Morgan have observed that the
inequality $\xi'(u)\ge \rho\cos\theta_v - \rho\cos\theta_u$
can be seen geometrically without a calculation.  This
inequality is all that is needed for the proof of Proposition
6.1-A.  (In fact, $\xi'\ge0$ is all that is needed.)
The heights of the two bumps in Diagram A-b are
$h(\theta_u)$, $h(\theta_v)$, where
$h(\theta) = \rho-\rho\cos\theta$.
We can increase the area by
$$(h(\theta_u)-h(\theta_v))\Delta u,$$
by cutting a vertical slice of area $h(\theta_v)\Delta u$
from the middle of the small bump, and adding a vertical slice of area
$h(\theta_u)\Delta u$ to the middle of the large bump.
To first order, this keeps the length of the perimeter constant.
The optimal increase in area $\Delta\xi$ is at least the increase
in area obtained by this strategy.  Hence the inequality.

The lemma leads to a proof of the Proposition.
By the equal curvature condition (6.2), if $u\ge v$,
then $\theta_u\ge\theta_v$, so lengthening longer chords decreases
the perimeter for fixed areas.

By continuity, it
is enough to prove the Proposition
when $\ell=a/b$ is a
rational number.
We apply the lemma to pairs of chords, increasing the longer
chord $u<s$ and decreasing the shorter chord $v>0$, keeping the sum
$u+v$ fixed,  continuing until every segment has length $0$ or $s$,
except for one of length between $0$ and $s$.

$\pi s^2/8$ is the area of a semicircle with diameter $s$.
By the equal curvature condition (6.2),
any circular arc greater
than a semicircle must lie along a chord of length $s$.
The area under such an arc is greater than $\pi s^2/8$, contrary
to hypothesis.

If $\sum \ell_i >\ell$, pick an edge of length $s$.  The
arc $e_i$ along that edge is less than a semicircle.
Decreasing the diameter $\ell_i$ while fixing $x(e_i)$ will decrease
the length of $e_i$.  (This is a standard argument: decrease
the obtuse angle
of the triangle joining the midpoint of $e_i$ to the endpoints of $f_i$.
This increases area of the triangle, keeping perimeter fixed.)
Continuing in this manner,
we can decrease $\sum \ell_i$ until
$\ell=\sum\ell_i$.  Again, we may assume that all lengths but one
are $s$ or $0$.

We replicate $b$ times these arcs and chords, enclosing
a total area of at most $b\pi s^2/8$.
We continue to apply the lemma to pairs of chords, until all edges
have length $0$ or $s$.
At no stage do the
circular arcs along the chords of length $s$ become semicircles.
The perimeter of the replicated version is $a\arc(s,b |X|/a)$,
so the perimeter of
the unreplicated version is $\ell \arc(s,|X(P)|/\ell)$
as desired.
\qed
\enddemo

\newpage
\head References
\endhead

\parskip=0.2\baselineskip

[A] F. J. Almgren, Jr.  Existence and regularity of almost everywhere
    of solutions to elliptic variational problems with constraints,
    Mem. AMS, 165 (1976).

[B] T. Bonnesen, Les probl\`emes des isop\'erim\`etres, 1929.

[BBC] A. Bezdek, K. Bezdek, R. Connelly, Finite and uniform stability
    of sphere packings, Discrete Comput Geom 20:111-130 (1998).

[C] C. Cox, L. Harrison, M. Hutchings, S. Kim, J. Light, A. Mauer, M.
    Tilton, The shortest enclosure of three connected areas in $\RR^2$,
    Real Analysis Exchange, Vol. 20(1), 1994/95, 313--335.

[CFG] H. Croft, K. Falconer, R. Guy, Unsolved Problems in Geometry,
    Springer, 1991.

[CS] J. H. Conway and N. J. A. Sloane,
    What are all the best sphere packings in low dimensions?
    Discrete Comput Geom 13:383-403 (1995).

[D] Charles Darwin, On the Origin of the Species.

[F] H. Federer, Geometric Measure Theory, Springer-Verlag, 1969.

[FT43] L. Fejes T\'oth, \"Uber das k\"urzeste Kurvennetz
    das eine Kugeloberfl\"ache in fl\"achengleiche konvexe
    Teil zerlegt, Mat. Term.-tud. \'Ertesit\"o 62 (1943), 349--354.

[FT64a] L. Fejes T\'oth, Regular Figures, MacMillan Company, 1964.

[FT64b] L. Fejes T\'oth, What the bees know and what they do not know,
    Bulletin AMS, Vol 70, 1964.

[He] T. Heath, A history of Greek mathematics, Vol II, Oxford, 1921.

[K] J. Kepler, L'\'etrenne ou la neige sexangulaire,
    Introduction by R. Halleux, C.N.R.S., 1975.

[Ku] G. Kuperberg, Notions of denseness, preprint,
    math.MG/9908003.

[L81] M. Lhuilier, M\'emoire sur le minimum de cire des
    alv\'eoles des abeilles, Nouveaux M\'emoires de l'Acad\'emie
    Royale des Sciences de Berlin, 1781.

[L89] M. Lhuilier, Abr\'eg\'e d'isop\'erim\'etrie \'el\'ementaire,
    1789.

[Mac] C. MacLaurin, Of the bases of the cells wherein the bees deposit
    their honey, Phil. Trans. Royal Society of London, 1743.

[M94] F. Morgan, Soap bubbles in ${\RR^2}$ and in surfaces,
    Pacific J. Math, 165 (1994), no. 2, 347--361.

[M95] F. Morgan, Geometric Measure Theory, A Beginner's Guide, Second
Edition, Academic Press, 1995.

[M99] F. Morgan, The hexagonal honeycomb conjecture,
    Trans. AMS, Vol 351, Number 5, pages 1753--1763,
    1999.

[P] Pappus d'Alexandrie,
    La collection math\'ematique,
    tr. Paul Ver Eecke, Albert Blanchard, 1982.

[Ph] R. Phelan, Generalisations of the Kelvin problem and other
    minimal problems, in [W], 1996.

[T] J. Taylor, The structure of singularities in soap-bubble-like
    and soap-film-like minimal surfaces, Annals of Math., 103 (1976),
    489-539.

[Th] D'Arcy Thompson, On Growth and Form, Cambridge, 1952.

[V] Marcus Terentius Varro, On Agriculture, Loeb Classical Library, 1934.

[W52] H. Weyl, Symmetry, Princeton, 1952.

[W] D. Weaire, The Kelvin problem: foam structures of minimal surface
    area, 1996.

[Wi] V. Willem, L'architecture des abeilles, 1928.

\bye